\documentclass[draft]{amsart}

\usepackage{amsmath}
\usepackage{amssymb}

\newcounter{bean}

\def\l{
\begin{list}
{\rm{(\alph{bean}).-}}{\usecounter{bean}
\setlength{\labelwidth}{0.8in} \setlength{\labelsep}{0.3cm}
\setlength{\leftmargin}{1cm}}}

\numberwithin{equation}{section}

\title[The Sequences of Fibonacci and Lucas for each Quadratic
Fields]{The Sequences of Fibonacci and Lucas for each real Quadratic
Fields $\mathbb{Q}(\sqrt{d}\ )$}

\author[Lam--Estrada]{Pablo Lam--Estrada}
\address{Escuela Superior de F\'isica y Matem\'aticas, Departamento de Matem\'aticas,
Instituto Polit\'ecnico Nacional (Unidad Zacatenco), CDMX, M\'exico}
\email{plam@esfm.ipn.mx}

\author[Maldonado--Ram\'irez]{Myriam Rosal\'ia Maldonado--Ram\'irez}
\address{Escuela Superior de F\'isica y Matem\'aticas, Departamento de Matem\'aticas,
Instituto Polit\'ecnico Nacional (Unidad Zacatenco), CDMX, M\'exico}
\email{rosalia@esfm.ipn.mx}

\author[L\'opez--Bonilla]{Jos\'e Luis L\'opez--Bonilla}
\address{Escuela Superior de Ingenier\'ia
Mec\'anica y El\'ectrica, Departamento de Ingenier\'ia en
Comunicaciones y Electr\'onica, Instituto Polit\'ecnico Nacional
(Unidad Zacatenco), CDMX, M\'exico} \email{jlopezb@ipn.mx}

\author[Jarqu\'in--Z\'arate]{Fausto Jarqu\'in--Z\'arate}
\address{Universidad Aut\'onoma de la Ciudad de M\'exico,
Academia de Matem\'aticas, Plantel San Lorenzo Tezonco, CDMX, M\'exico}
\email{fausto.jarquin@uacm.edu.mx}



\date{Abril 29, 2019}

\begin{document}

\begin{abstract}
We construct the sequences of Fibonacci and Lucas at any quadratic
field $\mathbb{Q}(\sqrt{d}\ )$ with $d>0$ square free, noting in
general that the properties remain valid as those given by the
classical sequences of Fibonacci and Lucas for the case $d = 5$,
under the respective variants. For this construction, we use the
fundamental unit of $\mathbb{Q}(\sqrt{d}\ )$ and then we observe the
generalizations for any unit of $\mathbb{Q}(\sqrt{d}\ )$ where,
under certain conditions, some of this constructions correspond to
$k$-Fibonacci sequence for some $k\in \mathbb{N}$. Of course, for
both sequences, we obtain the generating function, Golden ratio,
Binet's formula and some identities that they keep.
\end{abstract}

\maketitle

\section{Introduction}\label{S1}
The Fibonacci sequence was introduced by Leonardo of Pisa in 1202 in
his book Liber Abaci (Book of Calculation) [23]. Many of the
properties of the Fibonacci sequence were obtained by F. \'Edouard
Lucas who appoints such sequence by  ``Fibonacci'' [21, Section
$3{.}1{.}2$]. For more information about the history of the
Fibonacci numbers, we can see [20]. But also, Lucas is who initiates
the generalizations and their variants that have emerged from the
Fibonacci sequence, as we can observe, for example, in [4], [7],
[24], [25] and [26]. Vera W. de Spinadel introduced the Metallic
Means family whose members of such a family have many wonderful and
amazing properties, and applications to almost every areas of
sciences and arts, such as in some areas of the physical, biology,
astronomy and music (see [8], [9], [10] and [15]). On the other
hand, Sergio Falc\'on and \'Angel Plaza give properties of
$k$-Fibonacci sequence in [4], [5], [6] and [7], and these are a
particular case and general of metallic means families. Also in [3]
M. El-Mikkawy and T. Sogabe given a new family of $k$-Fibonacci
numbers.
 In [16], we can find hundreds of known identities, and
Azarian presents in [1] some known identities as binomial sums for quick numerical calculations.

In this paper, we associate with each real quadratic field
$\mathbb{Q}(\sqrt{d}\ )$, with $d>0$ square free, its owns sequences
of Fibonacci and Lucas (Definition 5), which correspond to certain
metallic means families (Theorem 8 and 13). These sequences of
Fibonacci and Lucas are determined by their generating functions
(Theorem 19) satisfying each Binet's formula (Theorem 22 and
Corollary 23). This means that each real quadratic field
$\mathbb{Q}(\sqrt{d}\ )$ will have also associated its own Golden
ratio (Definition 20), characteristic equation (\ref{eqn-11}) and
its Golden ratio will be the fundamental unit (Theorem 18). Finally,
we will establish for each $k\in \mathbb{N}$, the $k$-Fibonacci
sequence corresponds to Fibonacci sequence of the real quadratic
field $\mathbb{Q}(\sqrt{d}\ )$ for a unique $d>0$ square free
(Theorem 28).

At the time of submission, there is no description of the infinite
family of sequences in the Online Encyclopedia of Integer Sequences,
though some of the sequences do appear there, as indicated in Table
\ref{table1} and Table \ref{table2}.

This paper is organized as follows. In Section $2$ we collect
results of quadratic fields necessary for the development of the
work. In Section $3$ we construct the sequences of Fibonacci and
Lucas at any real quadratic field. Also we proof that the properties
remain valid as those given by the classical sequence of Fibonacci
and Lucas for $d=5$. In Section $4$ the main goal is proof that
Fibonacci and Lucas sequence are determined by the generating
functions. In Section $5$ we give Golden ratio associated as the
real quadratic field and we obtain Binet's formula in
$\mathbb{Q}(\sqrt{d}\;)$. In Section $6$ we extend our construct of
the sequences of Fibonacci and Lucas over all integer number.
Finally, in Section $7$ we define the sequence of Fibonacci and
Lucas of degree $d$ with respect to an arbitrary unit $\eta$ of
$\mathbb{Q}(\sqrt{d}\;)$ and we proof the results of the previous
sections are still met.

\section{Quadratic fields}
In this section we collect fundamental results from quadratic
fields. Throughout this paper, $d$ denotes a square free integer,
$\delta$ the discriminant of the quadratic field
$\mathbb{Q}(\sqrt{d}\ )$, $\mathcal{O}$ the ring of integers of
$\mathbb{Q}(\sqrt{d}\ )$, and $\mathcal{O}^*$ the multiplicative
group of all invertible elements of the ring $\mathcal{O}$. When
$d>0$, we say that $\mathbb{Q}(\sqrt{d}\ )$ is a {\bf real quadratic
field}, while if $d<0$ then $\mathbb{Q}(\sqrt{d}\ )$ is called an
{\bf imaginary quadratic field}. The following results are well
known.

\bigskip

\noindent {\bf Theorem 1.}\label{quadratic} {\it Keeping  the
previous notation.
\begin{enumerate}

\item[$(i)$] If $d\equiv 1$ {\rm mod} $4$, then the set $\displaystyle{\left \{\ 1\ ,\  \frac{1+\sqrt{d}}{2}\ \right\}}$
is an integral basis of $\mathbb{Q}(\sqrt{d}\ )$,\ $\delta=d$,\
$\mathcal{O}=\mathbb{Z}\displaystyle{\left[\frac{1+\sqrt{d}}{2}\right]}=
 \mathbb{Z}+\mathbb{Z}\ \displaystyle{\frac{1+\sqrt{d}}{2}}$\
 and\

 $\displaystyle{\mathcal{O}^* = \left \{ \ \frac{a+b\sqrt{d}}{2}\ \ \bigg |\ a,b\in \mathbb{Z},\ a^2-db^2=\pm 4 \right \} .}$

\item[$(ii)$] If $d\equiv 2$ {\rm mod} $4$ or
$d\equiv 3$ {\rm mod} $4$, then the set $\displaystyle{\left\{\ 1\ , \ \sqrt{d}\ \right\}}$
 is an integral basis of $\mathbb{Q}(\sqrt{d}\ )$,\  $\delta=4d$,\ $\displaystyle{\mathcal{O}=\mathbb{Z}[\sqrt{d}\ ]=
 \mathbb{Z}+\mathbb{Z}\ \sqrt{d}\ }$ and\

 $\displaystyle{\mathcal{O}^* =
 \left \{ \ a+b\sqrt{d}\ \ \big |\ a,b\in \mathbb{Z},\ a^2-db^2=\pm 1 \right \} .}$

\item[$(iii)$] If $d<0$, then $\mathcal{O}^* = \{\ -1, 1\ \}$
when $d\neq -1,\ -3$,\  $\mathcal{O}^* = \langle i \rangle = \{-1,
1, i, -i\}$ when $d=-1$\ and\ $\mathcal{O}^* = \langle \zeta_6
\rangle$ if $d=-3$, where $\zeta_6$ is a primitive $6$-th root of
unity.

\item[$(iv)$] If $d>0$, then

\begin{enumerate}
\item[$(a)$] There exists a unit $\varepsilon>1$ in $\mathcal{O}$ such that $\mathcal{O}^* = \langle -1 \rangle\times \langle \varepsilon\rangle $.

\item[$(b)$] If $u>1$ is a unit of $\mathcal{O}$, then $u=a+b\sqrt{d}$ for some $a>0$, $b>0$ in $\mathbb{Q}$.

\item[$(c)$] If $N(\varepsilon)=1$, then $N(u)=1$ for all $u\in\mathcal{O}^*$.
\end{enumerate}
\end{enumerate}
}

\noindent {\bf Proof.} See [13]. \hfill $\square$

\bigskip

The unit $\varepsilon$ of $\mathcal{O}$ in the Theorem 1, $(iv)$, is
called the {\bf fundamental unit} of $\mathcal{O}$. Hence the unit
$\varepsilon$ of $\mathcal{O}$ completely determines the group
$\mathcal{O}^*$. For example, we have for $d=5$ ($d\equiv 1$ mod
$4$),  $\displaystyle{\frac{1+\sqrt{5}}{2}}$ is the fundamental unit
of $\mathbb{Q}(\sqrt{5})$. If $d=17$, then
$\displaystyle{\frac{8+2\sqrt{17}}{2}= 4+\sqrt{17}}$ is a
fundamental unit of $\mathbb{Q}(\sqrt{17})$. In general, if $d\equiv
1$ mod $4$, then
$\varepsilon=\displaystyle{\frac{a_0+b_0\sqrt{d}}{2}}$ where $a_0$
and $b_0$ are either both even or both odd. Of course, if $a_0$ and
$b_0$ are both even, then $\varepsilon\in \mathbb{Z}[\sqrt{d}\ ]$.

\bigskip

On the other hand, we denote by $\mathfrak{M}_{2\times 2}(\mathbb{Z})$ the set of all matrices $2\times 2$ with integer entries.
Let $GL_2(\mathbb{Q})$ be the multiplicative group of invertible $2\times 2$ matrices
with rational entries, which is called the {\bf general lineal group of degree $2$ over $\mathbb{Q}$}.
The subset of all matrices of $GL_2(\mathbb{Q})$ with determinant $1$ is a normal subgroup of $GL_2(\mathbb{Q})$
called the {\bf special lineal group of degree $2$ over $\mathbb{Q}$} and denoted by $SL_2(\mathbb{Q})$.

\bigskip

For each $\lambda\in \mathbb{Q}$, let $$G_\lambda=\left\{A\in GL_2(\mathbb{Q})\ \left| \, A=\left[\begin{array}{cc}
a &  b\lambda \\
b  & a
\end{array}
\right]\right\}\right.\ , \ \ \ \ L_\lambda=\{A\in G_\lambda\ | ,
\det(A)=\pm 1\}\ \ \ {\rm and}$$

$$T_d=\left\{A\in\mathfrak{M}_{2\times 2}(\mathbb{Z})\ \left| \, A=\left[\begin{array}{cc}
a &  bd \\
b  & a
\end{array}
\right]\right\}\right. .$$

\bigskip

We have the follows results whose proofs can be seen in [17].

\bigskip

\noindent {\bf Theorem 2.}\label{5} {\it Keeping the previous
notation we obtain

\begin{enumerate}

\item[$(i)$] $T_d$ is a commutative subring with identity of $\mathfrak{M}_{2\times 2}(\mathbb{Z})$.

\item[$(ii)$] If $T^*_d$ is the multiplicative group of units of $T_d$, then  $T^*_d=L_d\cap \mathfrak{M}_{2\times 2}(\mathbb{Z})$. In particular, $T^*_d$ is a subgroup of $L_d$.

\item[$(iii)$] The rings $T_d$ and $\mathbb{Z}[\sqrt{d}\,]$ are isomorphic under the correspondence

$$\left[\begin{array}{cc}
a &  bd \\
b  & a
\end{array}
\right]   \longmapsto   a+b\sqrt{d}\ .$$  In particular, $T_d$ is an integral domain.

\item[$(iv)$] The isomorphism in $(iii)$ induces an isomorphism between the multiplicative groups $T^*_d$ and $(\mathbb{Z}[\sqrt{d}\, ])^*$.

\item[$(v)$] $T_d/(T_d\cap SL_2(\mathbb{Q}))\cong \{-1,\ 1\}$. \hfill $\square$

\end{enumerate} }

\bigskip

\noindent {\bf Theorem 3.}\label{6}
{\it Let $Q_d$ be the set of all matrices of the form
$A=\left[\begin{array}{cc}
a &  bd \\
b  & a
\end{array}
\right]$
with $a,b\in \mathbb{Q}$.

\begin{enumerate}

\item[$(i)$] $Q_d$ is a field isomorphic $\mathbb{Q}(\sqrt{d}\,)$ under the correspondence
$\left[\begin{array}{cc}
a &  bd \\
b  & a
\end{array}
\right]   \longmapsto   a+b\sqrt{d}\ .$ This is, $Q_d$ is the field of quotients of $T_d$.

\item[$(ii)$] There exists a monomorphism of the multiplicative group $\mathbb{Q}(\sqrt{d}\,)^*$ in the group $GL_2(\mathbb{Q})$.

\item[$(iii)$] The group $GL_2(\mathbb{Q})$ contains the chain of subgroups
$Q_d^*\cap SL_2(\mathbb{Q})<L_m<G_d=Q_d^*<GL_2(\mathbb{Q}).$

\hfill $\square$

\end{enumerate}
}

\noindent {\bf Theorem 4.}\label{7}
{\it Let $A=\left[\begin{array}{cc}
a &  bd \\
b  & a
\end{array}
\right]\in Q_d$  where $a$, $b$ are two rational numbers. Then the powers of $A$, $A^n=\left[\begin{array}{cc}
a_n &  b_nd \\
b_n  & a_n
\end{array}
\right]$ with $n\in \mathbb{N}$, are given as follows:

\begin{eqnarray}\label{eqn-1}
a_n & = & \left\{ \begin{array}{lcc}

\displaystyle{\sum_{0\leq t\leq \frac{n}{2}}\binom{n}{2t} a^{2t}b^{n-2t} d^{\frac{n}{2}-t}} & & if\ n\  even\\

& & \\

\displaystyle{\sum_{0\leq t\leq \frac{n-1}{2}}\binom{n}{2t+1} a^{2t+1}b^{n-2t-1} d^{\frac{n-1}{2}-t}} & &  if\ n\  odd
\end{array}\right .
\end{eqnarray}
and

\begin{eqnarray}\label{eqn-2}
b_n=\left\{ \begin{array}{lcc}

\displaystyle{\sum_{0\leq t\leq \frac{n-2}{2}}\binom{n}{2t+1} a^{2t+1}b^{n-2t-1} d^{\frac{n-2}{2}-t}} & & if\ n\ even\\

& & \\

\displaystyle{\sum_{0\leq t\leq \frac{n-1}{2}}\binom{n}{2t} a^{2t}b^{n-2t} d^{\frac{n-1}{2}-t}} & &  if\ n\  odd.
\end{array}\right .
\end{eqnarray}
\hfill $\square$ }

\section{The sequences of Fibonacci and Lucas in
$\mathbb{Q}(\sqrt{d}\ )$}
 In this section, we construct the
sequences of Fibonacci and Lucas at any real quadratic field. We
proof that the properties remain valid as those given by the
classical sequence of Fibonacci and Lucas for $d=5$. Being $d>0$ a
square free integer and $\varepsilon$ the fundamental unit of
$\mathbb{Q}(\sqrt{d}\ )$, we will write $\varepsilon= a+b\sqrt{d}$
where $a,b\in \mathbb{Q}$ with its corresponding matrix
$A_{\varepsilon}=\left[\begin{array}{cc}
a &  bd \\
b  & a
\end{array}
\right]\ $
and the powers $n$-th of $A_{\varepsilon}$ by
$A_{\varepsilon}^n=\left[\begin{array}{cc}
a_n &  b_nd \\
b_n  & a_n
\end{array}
\right]\ $ where $a_n$ and $b_n$ are given as in the equations
(\ref{eqn-1}) and (\ref{eqn-2}) of Theorem 4. Also, $\Delta_{}$ will
be the determinant of $A_{\varepsilon}$, that is,
$\Delta_{}=a^2-b^2d= N(\varepsilon)=\pm 1$, where $N$ is the {\bf
norm function} of the square field $\mathbb{Q}(\sqrt{d}\ )$.

\bigskip

Keeping the previous notation, we have the follows:

\bigskip

\noindent {\bf Definition 5.}\label{8}
{\it The sequence of \textbf{Fibonacci} $($resp.  \textbf{Lucas}$)$
\textbf{of degree d with respect to the fundamental unit
 $\varepsilon$} $($or simply the sequence of \textbf{Fibonacci} $($resp. \textbf{Lucas}$)$,
 if there is no risk of confusion with respect to $d$ and to its fundamental unit $\varepsilon)$
 is the sequence $\{F_{\varepsilon,n}\}_{n\in \mathbb{N}}$ $($resp. $\{L_{\varepsilon,n}\}_{n\in \mathbb{N}})$
 of positive numbers given as follows:
\begin{eqnarray}\label{eqn-3}
F_{\varepsilon,n}: =  \frac{b_n}{b}\ \ \ \ \ \left(resp.\ \
L_{\varepsilon,n}:= \frac{a_n}{a}\right)\ \ \ \ \  (n\in \mathbb{N})
\end{eqnarray}

\noindent where the sequence $\{b_n\}_{n\in \mathbb{N}}$ $($resp. $\{a_n\}_{n\in \mathbb{N}})$
is given as in the equation $(\ref{eqn-2})$ $($resp. $(\ref{eqn-1}))$ of Theorem $4$.
}

\bigskip

According to the equation (\ref{eqn-3}) of the Definition 5, we have that $F_{\varepsilon,n}$ and $L_{\varepsilon,n}$ are given by the follows equations:

\bigskip

\begin{eqnarray}\label{eqn-4}
F_{\varepsilon,n}=\left\{ \begin{array}{lcc}

\displaystyle{\sum_{0\leq t\leq \frac{n-2}{2}}\binom{n}{2t+1} a^{2t+1}b^{n-2t-2} d^{\frac{n-2}{2}-t}} & &  if\ n\  even\\

& & \\

\displaystyle{\sum_{0\leq t\leq \frac{n-1}{2}}\binom{n}{2t} a^{2t}b^{n-2t-1} d^{\frac{n-1}{2}-t}} & & if\ n\  odd
\end{array}\right .
\end{eqnarray}

\noindent and

\begin{eqnarray}\label{eqn-5}
L_{\varepsilon,n}=\left\{ \begin{array}{lcc}

\displaystyle{\sum_{0\leq t\leq \frac{n}{2}}\binom{n}{2t} a^{2t-1}b^{n-2t} d^{\frac{n}{2}-t}} & &  if\ n\ even\\

& & \\

\displaystyle{\sum_{0\leq t\leq \frac{n-1}{2}}\binom{n}{2t+1} a^{2t}b^{n-2t-1} d^{\frac{n-1}{2}-t}} & &  if\ n\ odd
\end{array}\right .
\end{eqnarray}

\noindent for each $n\in \mathbb{N}$.

\bigskip

In the Table \ref{table1} expresses some terms of the sequences
$\{F_{\varepsilon,n}\}_{n\in \mathbb{N}}$ and
$\{L_{\varepsilon,n}\}_{n\in \mathbb{N}}$ for some $d$'s square
free. Unless otherwise noted, the sequences are not in the Online
Encyclopedia of Integer Sequences at the time of publication, though
some of the sequences do appear there, as indicated in Table
\ref{table2}.

\newpage
\bigskip

\small{
\begin{table}[h!]
\centering
\begin{tabular}{|c|c|c|l|l|}
\hline
& & & & \\
$d$ & $\varepsilon$ & $\Delta$ & \hspace{0.25cm} Terms of the sequence & \hspace{0.25cm} Terms of the sequence  \\
& & & $ \hspace{1.4cm} \{F_{\varepsilon,n}\}_{n\in \mathbb{N}}$ & \hspace{1.4cm} $\{L_{\varepsilon,n}\}_{n\in \mathbb{N}}$ \\
& & & & \\
\hline
\hline
& & & &  \\
 $2$ & $1+\sqrt{2}$ & $-1$ & $1,2,5,12,29,70, 169, 408\ldots $ & $1,3,7,17,41,99,\ldots $ \\
 & &  &OEIS sequence A$000129$ & OEIS sequence A$001333$ \\
\hline
& & & & \\
$3$ & $2+\sqrt{3}$ & $1$ & $1,4,15,56,209,780,\ldots$ & $1,\displaystyle{\frac{7}{2}},13,\displaystyle{\frac{97}{2}},181,\displaystyle{\frac{1351}{2}},\ldots $  \\
& & & OEIS sequence A$001353$&OEIS sequence A$001075$ \\
\hline
& & & & \\
$5$ & $\displaystyle{\frac{1+\sqrt{5}}{2}}$ & $-1$ & $1,1,2,3,5,8,\ldots$ & $1,3,4,7,11,18,\ldots $  \\
& & & OEIS sequence A$000045$ & OEIS sequence A$000032$ \\
\hline
& & & &  \\
$6$ & $5+2\sqrt{6}$ & $1$ & $1, 10, 99, 980, 9701, 96030,\ldots$ & $1,\displaystyle{\frac{49}{5}},97,\displaystyle{\frac{4801}{5}},9505,\displaystyle{\frac{470449}{5}},\ldots $  \\
& & & OEIS sequence A$004189$& OEIS sequence A$001079$\\
\hline
& & & & \\
$7$ & $8+3\sqrt{7}$ & $1$ & $1,16,255,4064,64769,
\ldots$ & $1,\displaystyle{\frac{127}{8}},253,\displaystyle{\frac{32257}{8}},64261,
\ldots $   \\
& & & OEIS sequence A$077412$&OEIS sequence A$001081$ \\
\hline
& & & & \\
$10$ & $3+1\sqrt{10}$ & $-1$ & $1,6,37,228,1405,8658,\ldots$ & $1,\displaystyle{\frac{19}{3}},39,\displaystyle{\frac{721}{3}},1481,
\ldots $  \\
& & & OEIS sequence A$005668$& OEIS sequence A$005667$ \\
\hline
& & & & \\
$11$ & $10+3\sqrt{11}$ & $1$ & $1,20,399,7960,158801,
\ldots$ & $1,\displaystyle{\frac{199}{10}},397,\displaystyle{\frac{79201}{10}},158005,
\ldots $  \\
& &  & OEIS sequence A$075843$& OEIS sequence A$001085$ \\
\hline
& & & & \\
$13$ & $\displaystyle{\frac{3+\sqrt{13}}{2}}$ & $-1$ & $1,3,10,33,109,360,\ldots$ & $1,\displaystyle{\frac{11}{3}},12,\displaystyle{\frac{119}{3}},131,\displaystyle{\frac{1298}{3}},\ldots $  \\
& & &OEIS sequence A$006190$ & OEIS sequence A$006497$ \\
\hline
\end{tabular}
\caption{Sequence of Fibonacci and Lucas of degree $d$.}
\label{table1}
\end{table}
}

\small{
\begin{table}[h!]
\centering
\begin{tabular}{|c|c|c|l|c|}
\hline
& & & & \\
$d$ & $\varepsilon$ & $\Delta$ & \hspace{0.8cm} Terms of the sequence
 &
\hspace{0.08cm} OEIS   \\
& & & \hspace{1.8cm} $\{F_{\varepsilon,n}\}_{n\in \mathbb{N}}$ & integer sequence \\
& & & & \\
\hline \hline
& & & &  \\
 $37$ & $6+\sqrt{37}$ & $-1$ & $1,12,145,1752,21169,255780,\ldots $ & A$041061$ \\
 & &  &  &  \\
\hline
& & & &  \\
 $38$ & $37+6\sqrt{38}$ & $1$ & $1,74,5475,405076,29970149,2217385950,\ldots $ &  \\
 & &  &  &  \\
\hline
& & & &  \\
 $39$ & $25+4\sqrt{39}$ & $1$ & $1,50,2499,124900, 6242501, 312000150,\ldots $ &  \\
 & &  &  &  \\
\hline
& & & &  \\
 $41$ & $32 + 5\sqrt{41}$ & $-1$ & $1,64,4097,262272,16789505,1074790592,\ldots $ &  \\
 & &  &  &  \\
\hline
& & & &  \\
 $42$ & $13+2\sqrt{42}$ & $1$ & $1,26,675,17524, 454949, 11811150,\ldots $ & A$097309$ \\
 & &  &  &  \\
\hline
\end{tabular}
\caption{Sequence of Fibonacci of degree $d$.}
\label{table2}
\end{table}
}

\bigskip

\noindent {\bf Observation 6.}
Note that when $d=5$, we have   $\{F_{\varepsilon,n}\}_{n\in \mathbb{N}}$ and $\{L_{\varepsilon,n}\}_{n\in \mathbb{N}}$
are exactly the classical sequences of Fibonacci and Lucas, respectively.

\bigskip
\newpage

In the rest of the work, by abuse of notation, we write $F_n$ and
$L_n$ instead of $F_{\varepsilon,n}$ and $L_{\varepsilon,n}$ if
there is no risk of confusion with respect to the classical
sequences of Fibonacci and Lucas.

\bigskip

\noindent {\bf Theorem 7.}\label{t1}
{\it For each $m,n\in \mathbb{N}$,

\begin{enumerate}

\item[$(i)$] $F_{n+1}=a(L_{n}+F_{n})$.

\item[$(ii)$] $L_{n+1}=aL_{n}+\displaystyle{\frac{b^2d}{a}F_{n}}$.

\item[$(iii)$]

$F_{n}=\displaystyle{\frac{a}{\Delta}\left ( F_{n+1}-\frac{}{} L_{n+1} \right ) }=\left \{ \begin{array}{ccl}
a(L_{n+1}-F_{n+1}) & if & \Delta  =-1 \\
& & \\
a(F_{n+1}-L_{n+1}) & if & \Delta =1\ .
\end{array}
\right.$

\item[$(iv)$]

$L_{n}=\displaystyle{\frac{1}{\Delta}\left (  aL_{n+1}-\frac{b^2d}{a}F_{n+1}  \right )}=\left \{ \begin{array}{ccl}
\displaystyle{\frac{b^2d}{a}F_{n+1}-aL_{n+1}} & if & \Delta=-1 \\
& & \\
\displaystyle{aL_{n+1}-\frac{b^2d}{a}F_{n+1}}  & if & \Delta=1\ .
\end{array}
\right.$

\item[$(v)$] $F_{n+1}-a^n F_{1}= \displaystyle{\sum_{t=0}^{n-1}a^{t+1}L_{n-t}}$\ .

\item[$(vi)$] $L_{n+1}-a^n L_{1}= \displaystyle{b^2d\ \sum_{t=0}^{n-1}a^{t-1}F_{n-t}}$\ .

\item[$(vii)$] $F_{m+n}= a(F_{m} L_{n} + F_{n} L_{m})$.

\item[$(viii)$] $L_{m+n}= \displaystyle{\frac{b^2d}{a}}\cdot F_{m} F_{n} + a\ L_{m} L_{n} $.

\item[$(ix)$] $b^2d\ F_{n}^2-a^2L_{n}^2=-\Delta^n$.

\item[$(x)$] $F_n =\displaystyle{ \sum_{t=0}^{\lfloor \frac{n-1}{2} \rfloor} \binom {n}{2t+1} a^{n-2t-1} b^{2t} d^t} = \displaystyle{ \sum_{t=0}^{\lfloor \frac{n-1}{2} \rfloor} \binom {n}{n-2t-1} a^{n-2t-1} b^{2t} d^t}$.

\item[$(xi)$] $L_n =\displaystyle{ \sum_{t=0}^{\lfloor \frac{n}{2} \rfloor} \binom {n}{2t} a^{n-2t-1} b^{2t} d^t} = \displaystyle{ \sum_{t=0}^{\lfloor \frac{n}{2} \rfloor} \binom {n}{n-2t} a^{n-2t-1} b^{2t} d^t}$.

  \end{enumerate}

\noindent    Here $\lfloor x\rfloor$ is the integral part of $x\in \mathbb{R}$, i.e., is the greatest integer $n$ such that $n\leq x < n+1$.

\bigskip

}

\noindent {\bf Proof.} $(i)$ and $(ii)$ are obtained directly from
the equations (\ref{eqn-4}) and (\ref{eqn-5}). $(iii)$ and $(iv)$
are deducted from $(i)$ and $(ii)$. For induction, we obtain $(v)$
and $(vi)$. $(vii)$ and $(viii)$ are obtained from the relationship
$A_{\varepsilon}^{m+n}=A_{\varepsilon}^{m}\cdot
A_{\varepsilon}^{n}$. The relation
$\det(A_{\varepsilon}^n)=\Delta^n$ implies the relation $(ix)$.
Finally, $(x)$ y $(xi)$ are obtained of the relationships
$$a_n+b_n\sqrt{d}= (a+b\sqrt{d}\ )^n = \displaystyle{ \sum_{i=0}^n
\binom {n}{i} a^i b^{n-i} (\sqrt{d}\ )^{n-i}} = \displaystyle{
\sum_{i=0}^n \binom {n}{i} a^{n-i} b^{i} (\sqrt{d}\ )^{i}}.$$ \hfill
$\square$

\bigskip

\noindent {\bf Theorem 8.}\label{t2}
{\it There exist unique $r,s\in \mathbb{Q}^*$ such that $F_{n+2}=rF_{n}+sF_{n+1}$ for all $n\in \mathbb{N}$. More precisely, $F_{n+2}=(-\Delta) F_{n}+2a F_{n+1}$ for all $n\in \mathbb{N}$.}

\bigskip

\noindent {\bf Proof.} We have for each $n\in \mathbb{N}$,
\begin{eqnarray*}
(-\Delta) F_{n}+2a F_{n+1} & = & -(a^2-b^2d) F_n + 2a F_{n+1} = \displaystyle{a\left(\frac{b^2d}{a}\right ) F_n -a^2 F_n + 2a F_{n+1} }\\ & = & a(L_{n+1}-aL_n)-a^2 F_n +2a F_{n+1} = a(L_{n+1}+F_{n+1})\\ &=& F_{n+2}\ .
\end{eqnarray*}
On the other hand, let $r,s\in \mathbb{Q}^*$ be such that
\begin{eqnarray}\label{eqn-6}
F_{n+2}=rF_{n}+sF_{n+1}\ \ \ {\rm for\ all}\ n\in \mathbb{N}.
\end{eqnarray}
As $b^2d=a^2 - \Delta$, implies that $F_1=1$, $F_2= 2a$, $F_3= 4a^2-
\Delta$ and $F_4= 8a^3 - 4a  \Delta$. In particular, by the equation
(\ref{eqn-6}) for $n=1$ and $n=2$, we obtain the system of equations
\begin{eqnarray}\label{eqn-7}
\begin{array}{rcl}
r + 2as  & = & 4a^2- \Delta \\
& & \\

2ar + (4a^2 - \Delta) s & = &  8a^3 - 4a\Delta
\end{array}
\end{eqnarray}

\noindent which it has an unique solution, namely $r=- \Delta$ and $s=2a$. This complete the proof of theorem.

\hfill $\square$

\bigskip

\noindent {\bf Corollary 9.}\label{c1}
{\it The Fibonacci sequence $\{F_n\}_{n\in \mathbb{N}}$ is a $k$-Fibonacci sequence for some $k\in \mathbb{N}$ $($namely, $k=2a)$ if and only if $\Delta=-1$.}

\bigskip

\noindent {\bf Proof.} It immediate by Theorem 8. \hfill $\square$

\bigskip

\noindent {\bf Corollary 10.}\label{c2}
{\it The following conditions are equivalent:
\begin{enumerate}

\item[$(i)$] $F_{n+2}=F_{n}+F_{n+1}$ for all $n\in \mathbb{N}$;

\item[$(ii)$] $F_3=F_1+F_2$;

\item[$(iii)$] $d=5$ and $\displaystyle{\varepsilon=\frac{1+\sqrt{5}}{2}}$\ .

\end{enumerate}
}

\bigskip

\noindent {\bf Proof.}  $(i)\Longrightarrow (ii)$: It is immediate.

$(ii)\Longrightarrow (iii)$: We have that $-\Delta + 4a^2= (-\Delta)F_1+2a F_2 =F_3=F_1+F_2=1+2a$, then $4a^2-2a-(\Delta+1)=0$. If $\Delta=1$, then $2a^2-a-1=0$; since $a\neq 1$, necessarily $a=-1/2$. But this implies that $4b^2d=-3$; contradiction. Therefore $\Delta=-1$,  $a=1/2=b$ and $d=5$.

$(iii)\Longrightarrow (i)$: It is clear.
 \hfill $\square$

\bigskip

\noindent {\bf Corollary 11.}\label{c3}
{\it If $\{F_n\}_{n\in \mathbb{N}}$ is the  Fibonacci sequence classical, that is $d=5$, then
$$F_{n+2}=F_n+F_{n+1}$$
for each $n\in \mathbb{N}$.}

\bigskip

\noindent {\bf Proof.} It is immediate. \hfill $\square$

\bigskip

We recall if $d\equiv 2$ or $3$ mod $4$, then
$\varepsilon=a+b\sqrt{d}$ where $a,b\in \mathbb{Z}$. In this case,
it is obvious that $F_n\in \mathbb{N}$ for all $n\in \mathbb{N}$. If
$d\equiv 1$ mod $4$, then $\varepsilon =
a+b\sqrt{d}=\displaystyle{\frac{a_0+b_0\sqrt{d}}{2}}$ with $a_0,\
b_0\in \mathbb{N}$, where either are both even or both odd. When
they are both even, we have that $a, b\in \mathbb{N}$ and, hence,
$F_n\in \mathbb{N}$. But, in any case, $2a\in \mathbb{N}$.
Therefore, we obtain the following result.

\bigskip

\noindent {\bf Corollary 12.}\label{c41}
{\it $F_n\in \mathbb{N}$ for all $n\in \mathbb{N}$.}

\bigskip

\noindent {\bf Proof.} By Theorem 8, we have $F_{n+2}=(-\Delta)
F_{n}+2a F_{n+1}$ for all $n\in \mathbb{N}$, where $F_1=1$ and
$F_2=2a\in \mathbb{N}$. Then, the show follows by induction on $n$.
\hfill $\square$

\bigskip

\noindent {\bf Theorem 13.}\label{t3}
{\it There exist unique $r,s\in \mathbb{Q}^*$ such that $L_{n+2}=rL_{n}+sL_{n+1}$
for all $n\in \mathbb{N}$. More precisely, $L_{n+2}= (-\Delta)L_{n}+2a L_{n+1}$
 for all $n\in \mathbb{N}$.}

\bigskip

\noindent {\bf Proof.} We have  that for each $n\in \mathbb{N}$
\begin{eqnarray*}
(-\Delta) L_n + 2a L_{n+1} & = & \displaystyle{-\left ( a L_{n+1} - \frac{b^2d}{a} F_{n+1} \right )
+ 2a L_{n+1} } =  \displaystyle{ a L_{n+1} + \frac{b^2d}{a} F_{n+1}  }\\
& = & L_{n+2}\ .
\end{eqnarray*}
Now we prove the uniqueness. As $b^2d=a^2 - \Delta$, it follows that
$$\begin{array}{ccl}
L_1 & = & 1\\
L_2 & = & 2a - \displaystyle{\frac{\Delta}{a}}\\
L_3 & = & 4a^2 - 3\Delta\\
L_4 & = & 8a^3 - 8a\Delta + \displaystyle{\frac{1}{a}}\\
\vdots &  & \quad\quad\quad \vdots
\end{array}$$
Let $r,s\in \mathbb{Q}^*$ be such that
\begin{eqnarray}\label{eqn-8}
L_{n+2}=rL_{n}+sL_{n+1} \ \ \ {\rm for\ all}\ n\in \mathbb{N}.
\end{eqnarray}
In particular, for $n=1$ and $n=2$, we have the system of equations
\begin{eqnarray}\label{eqn-9}
\begin{array}{rcl}
r + \left (2a - \displaystyle{\frac{\Delta}{a}} \right ) s  & = & 4a^2- 3\Delta \\
& & \\

\left (2a - \displaystyle{\frac{\Delta}{a}} \right ) r + (4a^2 - 3\Delta) s & = &  8a^3 - 8a\Delta + \displaystyle{\frac{1}{a}}
\end{array}
\end{eqnarray}

\noindent which it has a unique solution, namely $r=-\Delta$ and
$s=2a$; so that, this system of equations has the same solution that
the system of equations (\ref{eqn-7}) given in the proof of Theorem
8. Therefore, the theorem is true. \hfill $\square$

\bigskip

Similarly to the corollaries to Theorem 8 for  Fibonacci sequence, we obtain corollaries to Theorem 13 for Lucas sequence.

\bigskip

\noindent {\bf Corollary 14.}\label{c4}
{\it The Lucas sequence $\{L_n\}_{n\in \mathbb{N}}$ is a $k$-Lucas sequence for some $k\in \mathbb{N}$ $($namely, $k=2a)$ if and only if $\Delta=-1$.}

\bigskip

\noindent {\bf Proof.} It immediate by Theorem 13. \hfill $\square$

\bigskip

\noindent {\bf Corollary 15.}\label{c5}
{\it The following conditions are equivalent:
\begin{enumerate}

\item[$(i)$] $L_{n+2}=L_{n}+L_{n+1}$ for all $n\in \mathbb{N}$;

\item[$(ii)$] $L_3=L_1+L_2$;

\item[$(iii)$] $d=5$ and $\displaystyle{\varepsilon=\frac{1+\sqrt{5}}{2}}$.

\end{enumerate}
}

\noindent {\bf Proof.}  $(i)\Longrightarrow (ii)$: It is immediate.

$(ii)\Longrightarrow (iii)$: Since $L_3=L_1+L_2$, that is, $4a^2-3\Delta=1+2a-\displaystyle{\frac{\Delta}{a}}$, we have that $4a^3-2a^2-a+(1-3a)\Delta=0$. If $\Delta=1$, then  $4a^3-2a^2-4a+1=0$ and $a$ can not be a rational number, contradiction. Hence, $\Delta=-1$ and $(2a^2+1)(2a-1)=0$. This implies that $a=1/2$ and $4b^2d=5$. Therefore, $d=5$ and $a=1/2=b$.

$(iii)\Longrightarrow (i)$: It is clear.
 \hfill $\square$

\bigskip

\noindent {\bf Corollary 16.}\label{c6}
{\it If $\{L_n\}_{n\in \mathbb{N}}$ is the Lucas sequence classical, that is $d=5$, then

$$L_{n+2}=L_n+L_{n+1}$$
for each $n\in \mathbb{N}$.}

\bigskip

\noindent {\bf Proof.} It is immediate. \hfill $\square$

\bigskip

\noindent {\bf Corollary 17.}\label{c61}
{\it For all $k\in \mathbb{N}$,

\begin{enumerate}
\item[$(i)$] $L_{2k-1}\in \mathbb{N}$;

\item[$(ii)$] if $a\in \mathbb{N}$, then $aL_{2k}\in \mathbb{N}$ and $(a,aL_{2k})=1$;

\item[$(iii)$] if $\displaystyle{a=\frac{a_0}{2}}$, with $a_0$ odd, then $a_0 L_{2k}\in \mathbb{N}$ and $(a_0,a_0L_{2k})=1$.

\end{enumerate}
}

\noindent {\bf Proof.} Applying the Theorem 13, the proof follows by
induction over all the pairs $(L_{2k-1},L_{2k})$, $k\in \mathbb{N}$.
\hfill $\square$

\section{Generating function}
The main goal of this section is to show that the Fibonacci and
Lucas sequences given in $(4)$ and $(5)$ are determined by the
generating functions.
\bigskip

\noindent {\bf Theorem 18.}\label{t4} {\it We obtain
\begin{enumerate}
\item[$(i)$] $\displaystyle{\lim_{n\rightarrow \infty}\frac{F_{n+1}}{F_n}}=\varepsilon= \displaystyle{\lim_{n\rightarrow \infty}\frac{L_{n+1}}{L_n}}$.

\item[$(ii)$] The series $\displaystyle{\sum_{n=1}^{\infty} F_n x^{n-1}}$ and $\displaystyle{\sum_{n=1}^{\infty} L_n x^{n-1}}$
both have the same radius of convergence, namely $R=1/\varepsilon$.
\end{enumerate}
}

\noindent {\bf Proof.} $(i)$: By Theorem 7, we have
$\displaystyle{\frac{F_{n+1}}{F_n}}=
\displaystyle{\frac{a(L_n+F_n)}{F_n}}= a+a\cdot
\displaystyle{\frac{L_n}{F_n}}= a+b\cdot
\displaystyle{\frac{a_n}{b_n}} $, and
$\displaystyle{\frac{L_{n+1}}{L_n}}=
\displaystyle{\frac{aL_n+\displaystyle{\frac{b^2d}{a}F_n}}{L_n}}=
a+\displaystyle{\frac{b^2d}{a}}\cdot \displaystyle{\frac{F_n}{L_n}}=
a+b\cdot \displaystyle{\frac{b_n d}{a_n}} $ where
$\displaystyle{\lim_{n\rightarrow \infty}\frac{a_n}{b_n}}=\sqrt{d}=
\displaystyle{\lim_{n\rightarrow \infty}\frac{b_{n} d}{a_n}}$, see
[17, Theorem 3.1]. Thus, $\displaystyle{\lim_{n\rightarrow
\infty}\frac{F_{n+1}}{F_n}}=\varepsilon=
\displaystyle{\lim_{n\rightarrow \infty}\frac{L_{n+1}}{L_n}}$.

$(ii)$: For each $x\in \mathbb{R}$, $x\neq 0$, we have that
$\displaystyle{\lim_{n\rightarrow\infty}\frac{F_{n+1}x^{n}}{F_{n}x^{n-1}}}=\varepsilon
x=
\displaystyle{\lim_{n\rightarrow\infty}\frac{L_{n+1}x^{n}}{L_{n}x^{n-1}}\
.}$ Then
$\displaystyle{\lim_{n\rightarrow\infty}\frac{F_{n+1}|x|^{n}}{F_{n}|x|^{n-1}}<1}$
if and  only if $\displaystyle{|x|<\frac{1}{\varepsilon}}$.
Similarly,
$\displaystyle{\lim_{n\rightarrow\infty}\frac{L_{n+1}|x|^{n}}{L_{n}|x|^{n-1}}<1}$
if and  only if $\displaystyle{|x|<\frac{1}{\varepsilon}}$\ .
Therefore, both series have the same radius of convergence
$R=1/\varepsilon$. This complete the proof of the theorem. \hfill
$\square$

\bigskip

\noindent {\bf Theorem 19.}\label{t5} {\bf (Generating function)}
{\it Let $x\in \mathbb{R}$ be such that $|x|<1/\varepsilon$.
\begin{enumerate}
\item[$(i)$]  If $f(x)=\displaystyle{\sum_{n=1}^{\infty} F_n x^{n-1}}$, then $\displaystyle{f(x)=\frac{1}{\Delta x^2-2ax+1}\ .}$

\item[$(ii)$] If $g(x)=\displaystyle{\sum_{n=1}^{\infty} L_n x^{n-1}}$, then $\displaystyle{g(x)=\left ( \frac{a-\Delta x}{a}\right ) f(x)= \frac{a-\Delta x}{a(\Delta x^2-2ax+1)}\ .}$
\end{enumerate}
}

\noindent {\bf Proof.} $(i)$: For each $x\in \mathbb{R}$ with $|x|<1/\varepsilon$, we have that

\begin{eqnarray*}
f(x) & = & \sum_{n=1}^{\infty} F_n x^{n-1} = 1+ 2a x + \sum_{n=1}^{\infty} F_{n+2} x^{n+1}  =  \displaystyle{1+ 2ax + \sum_{n=1}^{\infty} \left ( - \Delta F_n + 2a F_{n+1} \right ) x^{n+1}}\\
& = & 1+2ax  - \Delta x^2 f(x) + 2ax ( f(x)-1)  =  1+ f(x) ( 2ax -\Delta x^2 )
\end{eqnarray*}

\noindent this implies that
$\displaystyle{f(x)=\frac{1}{\Delta x^2 -2ax +1}\ .}$

$(ii)$: We observe that, for each $x\in \mathbb{R}$ with $|x|<1/\varepsilon$
\begin{eqnarray*}
f(x) & = & 1 + \sum_{n=1}^{\infty} F_{n+1} x^{n}=1 + \sum_{n=1}^{\infty} a(L_n+F_n)x^n =  1+ ax g(x)+ ax f(x)\ .
\end{eqnarray*}
If $x\neq 0$, then
\begin{eqnarray*}
g(x) & = & \frac{(1-ax)f(x)-1}{ax} = \frac{a-\Delta x}{a(\Delta x^2 -2ax +1)} =  \left( \frac{a-\Delta x}{a} \right ) f(x)\ .
\end{eqnarray*} \hfill $\square$

\section{Golden ratio and Binet's formula in $\mathbb{Q}(\sqrt{d}\
)$}

In this section we give Golden ratio associated as the quadratic
field $\mathbb{Q}(\sqrt{d}\,)$. Also we obtain Binet's formula in
$\mathbb{Q}(\sqrt{d}\,)$. We start with

\bigskip

\noindent {\bf Definition 20.}\label{ratio} {\it Let $x, y\in
\mathbb{R}$ be such that $0<y<x$. We say that $x$ and $y$ are in
\textbf{Golden ratio with respect to the quadratic field
$\mathbb{Q}(\sqrt{d}\, )$} $($or simply that they  are in
\textbf{Golden ratio},  if there is no risk of confusion with
respect to the quadratic field $\mathbb{Q}(\sqrt{d}\ ))$,  if
\begin{eqnarray}\label{eqn-10}
\frac{2ax-\Delta y}{x}=\frac{x}{y}\ .
\end{eqnarray}
}

\bigskip

Thus, if $x$ and $y$ are in Golden ratio and we write $\varphi:=\displaystyle{\frac{x}{y}}$, then we have that

\begin{eqnarray*}
2a-\frac{\Delta}{\varphi} & =  & 2a-\Delta\cdot \frac{y}{x}=\frac{2ax-\Delta y}{x} =\frac{x}{y}\\ & = & \varphi.
\end{eqnarray*}
This is, $\varphi$ satisfies the equation
\begin{eqnarray}\label{eqn-11}
\varphi^2-2a  \varphi + \Delta =0.
\end{eqnarray}

But $x^2-2ax+\Delta$ is the irreducible polynomial of $\varepsilon$ over $\mathbb{Q}$ with $\overline{\varepsilon}$ its other root, where $\overline{\varepsilon}$ is the conjugate of $\varepsilon$. Therefore, $\varphi=\varepsilon$ or $\varphi=\overline{\varepsilon}$. As $x>y>0$ and $\overline{\varepsilon}=\Delta/\varepsilon$, necessarily $\varphi=\varepsilon$. In consequence, we have the equation
\begin{eqnarray}\label{eqn-12}
\varepsilon^2=2a\varepsilon-\Delta.
\end{eqnarray}

\bigskip

\noindent {\bf Theorem 21.}\label{t6}
{\it For each  $n\in \mathbb{N}$, with $n\geq 2$,
\begin{eqnarray}\label{eqn-13}
\varepsilon^{n}=F_{n}\ \varepsilon-F_{n-1}\ \Delta .
\end{eqnarray}
}

\noindent {\bf Proof.} The show is by induction on $n$. It is clear
for $n=2$, that is, $\varepsilon^{2}=2a \varepsilon-\Delta=F_{2}\
\varepsilon-F_{1}\ \Delta$. Hence,
\begin{eqnarray*}
\varepsilon^{n+1} & = & \varepsilon (F_n\ \varepsilon-F_{n-1}\
\Delta)=F_n(2a\varepsilon-\Delta)-F_{n-1} \varepsilon \Delta =
(-\Delta F_{n-1}+2aF_n)\varepsilon-F_n \Delta\\ &=& F_{n+1}\
\varepsilon - F_n \Delta.
\end{eqnarray*} \hfill $\square$

\bigskip

Since $\overline{\varepsilon}$ also satisfies the equation (\ref{eqn-12}), we have the equation
\begin{eqnarray}\label{eqn-14}
(\overline{\varepsilon})^{n}=F_{n}\ \overline{\varepsilon}-F_{n-1}\ \Delta,
\end{eqnarray}

\noindent for each $n\geq 2$.

\bigskip

\noindent {\bf Theorem 22.}\label{Binet}
{\it For each $n\in \mathbb{N}$,
\begin{eqnarray}\label{eqn-15}
F_{n}=\frac{\ \varepsilon^n - (\overline{\varepsilon})^n}{\varepsilon - \overline{\varepsilon}}\ .
\end{eqnarray}
}

\noindent {\bf Proof.} It follows to make the difference of the equations (\ref{eqn-13}) and (\ref{eqn-14}). \hfill $\square$

\bigskip

The equation (\ref{eqn-15}) is known as the {\bf Binet's formula}

\bigskip

\noindent {\bf Corollary 23.}\label{c7}
{\it For each $n\in \mathbb{N}$,
\begin{eqnarray}\label{eqn-16}
L_n=\frac{\varepsilon^n+(\overline{\varepsilon})^{n}}{\varepsilon+\overline{\varepsilon}}\ .
\end{eqnarray}
}

\noindent {\bf Proof.} It is immediate from the following

\begin{eqnarray*}\label{eqn-17}
\varepsilon^n+(\overline{\varepsilon})^{n} &= & 2aF_n-2\Delta F_{n-1}= 2a\left (F_n-\frac{\Delta}{a}\cdot F_{n-1}\right ) =
 2a L_n =(\varepsilon+\overline{\varepsilon})L_n.
\end{eqnarray*} \hfill $\square$

\bigskip

The following two theorems give us other version of the generating
functions of the sequences of  Fibonacci and   Lucas in
$\mathbb{Q}(\sqrt{d}\ )$.

\bigskip

\noindent {\bf Theorem 24.}\label{tfgf}
{\it Let $f_1(x)=\displaystyle{\sum_{n=0}^{\infty} \Delta^{n} F_{n+1}x^n}$ and $\displaystyle{g_1(x)=\sum_{n=0}^{\infty} \Delta^{n} L_{n+1}x^n}$. Then, the series $f_1(x)$ and $g_1(x)$ are convergent for $|x|<{\rm min}\{|\varepsilon|,|\overline{\varepsilon}|\}$. Furthermore,
\begin{eqnarray}\label{fgf}
f_1(x) & = & \frac{\Delta}{x^2-2a x+ \Delta}
\end{eqnarray}
\noindent and
\begin{eqnarray}\label{fgl}
g_1(x) & = & \frac{\Delta(a-x)}{a(x^2-2a x+ \Delta)} = \left ( \frac{a-x}{a}\right ) f_1(x).
\end{eqnarray}
}

\noindent {\bf Proof.} We have for $|x|<{\rm
min}\{|\varepsilon|,|\overline{\varepsilon}|\}$

\begin{eqnarray*}
\frac{2b\sqrt{d}}{(x-\varepsilon)(x-\overline{\varepsilon})} & = & \frac{1}{x-\varepsilon}-\frac{1}{x-\overline{\varepsilon}} = \frac{1}{\overline{\varepsilon}\left ( 1-\displaystyle{\frac{x}{\overline{\varepsilon}}}\right)} - \frac{1}{\varepsilon\left ( 1-\displaystyle{\frac{x}{\varepsilon}}\right)}\\
& = & \frac{1}{\overline{\varepsilon}}\ \sum_{n=0}^{\infty} \left ( \frac{x}{\overline{\varepsilon}} \right )^n - \frac{1}{\varepsilon}\ \sum_{n=0}^{\infty} \left ( \frac{x}{\varepsilon} \right )^n = \sum_{n=0}^{\infty} \left (\frac{\varepsilon^{n+1}-\overline{\varepsilon}^{\ n+1}}{(\varepsilon \overline{\varepsilon})^{n+1}}\right )x^n\\
& = & 2b\sqrt{d} \left (\sum_{n=0}^{\infty} (a^2-b^2d)^{n+1} \left (\frac{\varepsilon^{n+1}-\overline{\varepsilon}^{\ n+1}}{\varepsilon - \overline{\varepsilon}}\right )x^n \right)\\
& = & 2b\sqrt{d} \left (\sum_{n=0}^{\infty} \Delta^{n+1}F_{n+1}x^n \right).
\end{eqnarray*}
This implies that
\begin{eqnarray*}
\frac{1}{x^2-2a x+ \Delta} & = & \frac{1}{x^2-2a x+ (a^2-b^2d)}  = \frac{1}{(x-\varepsilon)(x-\overline{\varepsilon})}\\ & = & \sum_{n=0}^{\infty} \Delta^{n+1} F_{n+1}x^n
\end{eqnarray*}
 or equivalently
\begin{eqnarray}\label{eqn-18}
\frac{\Delta}{x^2-2a x+ \Delta} = \sum_{n=0}^{\infty} \Delta^{n+2} F_{n+1}x^n = \sum_{n=0}^{\infty} \Delta^{n} F_{n+1}x^n.
\end{eqnarray}
Therefore
\begin{eqnarray*}
f_1(x) & = & \frac{\Delta}{x^2-2a x+ \Delta}\ .
\end{eqnarray*}
On the other hand, we have
\begin{eqnarray*}
\left ( a -x\right )f_1(x) &=&  \left ( a -x\right )  \left (\sum_{n=0}^{\infty} \Delta^{n+2} F_{n+1}x^n\right )\\
&=& a\ F_1 + \sum_{n=1}^{\infty}  \Delta^{n+1} \left ( \frac{a}{\Delta}  F_{n+1}-  F_{n} \right ) x^n\\
&=& a\ L_1 + \sum_{n=1}^{\infty}  \Delta^{n} a  L_{n+1} x^n\\ & = & a\ \sum_{n=0}^{\infty}  \Delta^{n}   L_{n+1} x^n.
\end{eqnarray*}
Therefore
\begin{eqnarray*}
g_1(x) & = & \left ( \frac{a-x}{a}\right ) f_1(x)=\frac{\Delta(a-x)}{a(x^2-2a x+ \Delta)}\ .
\end{eqnarray*}
\hfill $\square$

\section{Some Other Properties}
Using the equations (\ref{eqn-15}) and (\ref{eqn-16}), we can extend
the definition of the sequences of Fibonacci and Lucas over all
integer number. This is, we use the Binet's formula for all $n\in
\mathbb{Z}$, Theorem 22 and Corallary 23 we obtain
\begin{eqnarray}
F_{-n} &= &\left \{ \begin{array}{ccl}
0 & if & n=0\\
& & \\
-\Delta ^n F_n & if & n\geq 1
\end{array} \right .
\end{eqnarray}

\noindent and

\begin{eqnarray}
L_{-n} &= &\left \{ \begin{array}{ccl}
\displaystyle{\frac{1}{a}} & if & n=0\\
& & \\
\Delta ^n L_n & if & n\geq 1.
\end{array} \right .
\end{eqnarray}
Thus, it holds that for all $n\in \mathbb{Z}$
\begin{eqnarray}
F_{n+2} & = & (-\Delta) F_n + 2a F_{n+1}
\end{eqnarray}
and
\begin{eqnarray}
L_{n+2} & = & (-\Delta) L_n + 2a L_{n+1}.
\end{eqnarray}

\bigskip

But also we can obtain, in our case, the identities established by
Catalan, Cassini, D' Ocagne, and Hosnberger which are hold for all
$n\in \mathbb{Z}$, that is

\bigskip

\noindent {\bf Theorem 25.} {\it For all $m,n\in \mathbb{Z}$, the follows identities are holds:

\begin{enumerate}

\item[$(i)$] $F_n^2-F_{n+m}F_{n-m}=\Delta^{n-m} F_{m}^2.$

\item[$(ii)$] $F_n^2-F_{n-1} F_{n+1}  = \Delta^{n-1}.$

\item[$(iii)$] $L_n^2-L_{n+r}L_{n-r}=\displaystyle{\frac{\Delta^n}{2a^2}-\left(\frac{\Delta^{n-r}}{2a}\right) L_{2r}} $.

\item[$(iv)$] $F_m F_{n+1} - F_n F_{m+1} =  \Delta^{n} F_{m-n}.$

\item[$(v)$] $F_{m-1} F_n + F_m F_{n+1}= \left \{ \begin{array}{lcl}
F_{m+n} &  if & \Delta = -1 \\
& & \\
\displaystyle{\frac{a}{2b^2d}}\cdot \bigg (2a L_{m+n} - L_{m-n-1}\bigg ) &  if & \Delta = 1.
\end{array}
\right. $

\item[$(vi)$] $L_n L_{n+r} = \displaystyle{\left (\frac{1}{2a}\right )\ L_{2n+r} + \left (\frac{\Delta^n}{2a}\right ) \ L_r}$.

\end{enumerate}
}

\noindent {\bf Proof.} The show for each of the identities can be
performed using the Binet's Formula. So we will prove only $(iv)$.
Hence we have
\begin{eqnarray*}
F_m F_{n+1} - F_n F_{m+1} & = & \left( \frac{\varepsilon^m-(\bar{\varepsilon})^m}{\varepsilon-\bar{\varepsilon}}\right)
\left( \frac{\varepsilon^{n+1} -(\bar{\varepsilon})^{n+1}}{\varepsilon-\bar{\varepsilon}}\right)\\
& & - \left( \frac{\varepsilon^{n}- (\bar{\varepsilon})^{n}}{\varepsilon-\bar{\varepsilon}}\right) \left( \frac{\varepsilon^{m+1}
 -(\bar{\varepsilon})^{m+1}}{\varepsilon-\bar{\varepsilon}}\right) \\
 & = & \frac{\varepsilon^m (\bar{\varepsilon})^{n} -
\varepsilon^{n}(\bar{\varepsilon})^{m}}{\varepsilon-\bar{\varepsilon}}
 =  \frac{\varepsilon^{m-n+n} (\bar{\varepsilon})^{n} -
\varepsilon^{n}(\bar{\varepsilon})^{m-n+n}}{\varepsilon-\bar{\varepsilon}} \\
& = & \Delta^{n} \left( \frac{\varepsilon^{m-n}  -
(\bar{\varepsilon})^{m-n}}{\varepsilon-\bar{\varepsilon}}\right)  =
\Delta^{n} F_{m-n}.
\end{eqnarray*}

 \hfill $\square$

\section{The sequence of Fibonacci and of Lucas of degree $d$ with
respect to an arbitrary unit}

The unit group of $\mathbb{Q}(\sqrt{d}\ )$, with $d>0$,  is isomorph to group $\langle -1 \rangle\times  \langle \varepsilon \rangle$ where $\varepsilon$
is the fundamental unit of $\mathbb{Q}(\sqrt{d}\ )$, generator of the  infinity cyclic subgroup.
This cyclic subgroup also is generated by $1/\varepsilon$, $-\varepsilon$ and $-1/\varepsilon$.
 This is, each unit of $\mathbb{Q}(\sqrt{d}\ )$ has the form $\pm \varepsilon^l$ for some $l\in \mathbb{Z}$.
 Observing the previous development, we can define the sequence of Fibonacci  and  Lucas of degree $d$ with respect
 to an arbitrary unit $\eta$ of $\mathbb{Q}(\sqrt{d}\ )$, and the results of the previous sections are still met.
 Essentially this is because $N(\eta)=\pm 1$. This allows us to build even more an infinity
 of sequences in $\mathbb{Q}(\sqrt{d}\ )$ meeting similar properties of the sequences of Fibonacci and Lucas.
  For example, we consider the unit $\eta=\displaystyle{\frac{-1+\sqrt{5}}{2}}$ of $\mathbb{Q}(\sqrt{5}\ )$,
  we have that the first terms of the sequence of Fibonacci of degree $5$ with respect to the unit $\eta$  are:

$$F_{\eta,1}=1,\ F_{\eta,2}=-1,\ F_{\eta,3}=2,\ F_{\eta,4}=-3,\ \ldots$$

\noindent where $N(\eta)=-1$. Comparing the terms of the sequence of Fibonacci with negative index,
$F_{-n}$ with $n\geq 1$, we have that $F_{\eta,n}=F_{-n}$ for all $n\in \mathbb{N}$. This is,
the sequence of Fibonacci with negative index
 of degree $5$ with respect to the fundamental unit $\displaystyle{\varepsilon=\frac{1+\sqrt{5}}{2}}$
 is the sequence of Fibonacci of degree $5$ respect to the unit $\eta=\displaystyle{\frac{-1+\sqrt{5}}{2}}$\ . This it is not
 a coincidence, that is, this fact is generalized in the following.

 \bigskip

\noindent {\bf Theorem 26.}
{\it The Fibonacci sequence of degree $d$ with respect to the unit $1/\varepsilon$
and $\Delta=-1$ is the Fibonacci sequence with negative index of degree $d$ with respect to the fundamental unit $\varepsilon$.}

\bigskip

\noindent {\bf Proof.} We write $\eta=1/\varepsilon=\Delta \overline{\varepsilon}$.
Hence,  $\overline{\eta}= \Delta\varepsilon$. Using the Binet's formula, we have that for all $n\in \mathbb{N}$
\begin{eqnarray*}
F_{\eta, n} & = & \frac{\eta^n-(\overline{\eta})^n}{\eta-\overline{\eta}}
=\frac{(\Delta\overline{\varepsilon})^n-(\Delta\varepsilon)^n}{\Delta\overline{\varepsilon}-\Delta\varepsilon}
 =  \Delta^{n-1}\ \frac{\varepsilon^n-(\overline{\varepsilon})^n}{\varepsilon-\overline{\varepsilon}}=\Delta^{n-1}F_n=-\Delta^n F_n
 =  F_{-n}\ .
\end{eqnarray*}
\hfill $\square$

\bigskip

\noindent {\bf Observation 27.}
We note that if $\Delta=1$ then the  Fibonacci  sequence of degree $d$ with respect to the unit $1/\varepsilon$
coincides with the Fibonacci  sequence of degree $d$ with respect to the unit $\varepsilon$.

\bigskip

We finish our work with the following result.

\bigskip

\noindent {\bf Theorem 28.}\label{k-1} {\it For each $k\in
\mathbb{N}$ there exist unique $d,\ r\in \mathbb{N}$ such that $d$
is square free and $\displaystyle{\frac{k}{2}+\frac{r}{2}\sqrt{d}}$
is a unit of the quadratic field $\mathbb{Q}(\sqrt{d}\, )$ with norm
$-1$. Therefore, in this case, the $k$-Fibonacci sequence is the
Fibonacci  sequence of degree $d$ with respect to a unit of
$\mathbb{Q}(\sqrt{d}\, )$.}

\bigskip

\noindent {\bf Proof.} Let $k\in \mathbb{N}$ be arbitrary.
We have that $k^2+4$ is not a perfect square.
Hence, there exist $d, r\in \mathbb{N}$ such that $k^2+4=r^2 d$ where $d$ is positive square free.
 This implies that $\displaystyle{\left (\frac{k}{2}\right )^2-\left (\frac{r}{2}\right )^2d=-1}$.
 Hence, $\displaystyle{\frac{k}{2}+\frac{r}{2}\sqrt{d}}$ is a unit of $\mathbb{Q}(\sqrt{d}\ )$
 with norm $-1$. On the other hand, if $d,\ d_1, r,\ r_1\in \mathbb{N}$
 such that $\displaystyle{\frac{k}{2}+\frac{r}{2}\sqrt{d}}$ and
  $\displaystyle{\frac{k}{2}+\frac{r_1}{2}\sqrt{d_1}}$ are units of the
  quadratic field $\mathbb{Q}(\sqrt{d}\ )$ both with norm $-1$,
  then $\displaystyle{\left (\frac{k}{2}\right )^2-\left (\frac{r}{2}\right )^2d=-1=\left
  (\frac{k}{2}\right )^2-\left (\frac{r_1}{2}\right )^2d_1 }$,
  thus $\displaystyle{\left (\frac{r}{2}\right )^2d=\left (\frac{r_1}{2}\right )^2d_1 }$.
  That is, $r^2d=r_1^2d_1$, where $d$ and $d_1$ are square free. Therefore, $d_1=d$ and $r_1=r$.
  In consequence, the $k$-Fibonacci sequence is the Fibonacci sequence of degree $d$
  with respect to a unit of $\mathbb{Q}(\sqrt{d}\ )$. \hfill $\square$

\bigskip

\noindent {\bf Corollary 29.} {\it For each $k\in \mathbb{N}$, the
$k$-Fibonacci sequence is the Fibonacci sequence of degree $d$ with
respect to a unit of $\mathbb{Q}(\sqrt{d}\ )$ for some $d$ square
free.}

\bigskip

\noindent {\bf Proof.} It is immediately of Theorem 28 and Corollary
9. \hfill $\square$

\section{Conclusions}
In this work we have established that every real quadratic field
$\mathbb{Q}(\sqrt{d}\ )$ has its own Fibonacci sequence and Lucas
sequence, and variants of these, through the fundamental unit, being
this the golden ratio. Therefore, the real quadratic field
$\mathbb{Q}(\sqrt{d}\ )$ has its own gold ratio. Under these
conditions, it is possible that may arise further research aimed at
obtaining properties, both algebraic and geometric, related with the
intrinsic properties of the real quadratic field
$\mathbb{Q}(\sqrt{d}\ )$.

\bigskip


\noindent {\bf References}

\bigskip

\begin{itemize}

\item [{[1]}] Azarian  M. K., {\it Identities Involving Lucas or Fibonacci and Lucas Numbers as Binomial Sums}.
 Int. J. Contemp. Math. Sciences, {\bf 7}, (2012) No. 45, 2221-2227.

\item [{[2]}]   Bolat C. and K\"ose H., {\it On the Properties of $k$-Fibonacci Numbers.} Int. J. Contemp. Math. Sciences, {\bf 5}, (2010) No. 22, 1097-1105.

\item [{[3]}]   El-Mikkawy M.  and Sogabe T., {\it A new family of $k$-Fibonacci numbers}. Applied Mathematics and Computation, {\bf 215}, (2010) 4456-4461.

\item [{[4]}]   Falc\'on S., {\it On the $k$-Lucas Numbers}. Int. J. Contemp Math. Sciences, {\bf 6}, (2011) No. 21, 1039-1050.

\item [{[5]}]   Falc\'on S. and Plaza A., {\it The $k$-Fibonacci sequence and the Pascal $2$-triangle}. Chaos, Solitons and Fractals, {\bf 33}, (2007) 38-49. 

\item [{[6]}]   Falc\'on S. and Plaza A., {\it The $k$-Fibonacci hyperbolic functions}. Chaos, Solitons and Fractals, {\bf 38}, (2008) 409-420. 

\item [{[7]}]   Falc\'on S. and Plaza A., {\it On $k$-Fibonacci sequences and polynomials and their derivatives}. Chaos, Solitons and Fractals, {\bf 39}, (2009) 1005-1019. 

\item [{[8]}]   Fujita M. and  Machida K., {\it Spectral properties of onedimensional
quasi-crystalline and incommensurate systems}, J.
Phys. Soc. Jpn, {\bf 56}, (1987), No. 4, 1470-1477. 

\item [{[9]}]    Gumbs G. and  Ali M. K., {\it Dynamical maps, Cantor spectra and
localization for Fibonacci and related quasiperiodic lattices}, Phys. Rev.
Lett., {\bf 60}, (1988), No. 11,  1081-1084. 

\item [{[10]}]    Gumbs G. and  Ali M. K., {\it Electronic properties of the tight-binding
Fibonacci Hamiltonian}, J. Phys. A: Math. Gen., {\bf 22}, (1989), No. 8,  951-970. 

\item [{[11]}]    Garc\'ia Hern\'andez V. C. and  Mej\'ia Huguet V. J., {\it N\'umeros de Fibonacci}. UAM, M\'exico, 2014. 

\item [{[12]}]    Hungerford T. W., {\em Algebra},
Graduate Texts in Mathematics {\bf 73}, Springer-Verlag New York, Inc., 1974. 

\item [{[13]}]    Janusz G. J., {\it Algebraic Number Fields}.
Graduate Studies in Mathematics, {\bf 7}. AMS, 1996. 

\item [{[14]}]    Jhala D.,  Rathore G.P.S. and  Singh B., {\it Some Identities Involving Common Factors of k-Fibonacci and k-Lucas
Numbers},
American Journal of Mathematical Analysis, {\bf 2}, (2014)  No. 3, 33-35. 

\item [{[15]}]    Kappraff J., {\it Musical proportions at the basis of systems of
architectural proportions}, NEXUS: Architecture and Mathematics, Edizioni dell'Erba, 1996. 

\item [{[16]}]   Koshy T., {\it Fibonacci and Lucas Numbers with Applications.} John Wiley \& Sons, Inc., 2001. 

\item [{[17]}]    M\'endez-Delgadillo H.,  Lam-Estrada P. and  Maldonado-Ram\'irez M. R., {\it Approach to Square Roots Applying Square Matrices.}
Palestine Journal of Mathematics, {\bf 4}, (2015), No. 2, 271-276. 

\item [{[18]}]    Mollin R. A., {\it Algebraic Number Theory}, Chapman and Hall/CRC Press. 1999. 

\item [{[19]}]    Niven I. M.,  Zukermann H. S. and  Montgomery H. L., {\em An introduction to the theory of numbers},
Wiley, 1991. 

\item [{[20]}]    Posamentier A. S. and  Lehmann I., {\it The fabulous Fibonacci numbers}. Amherst, NY: Prometheus Books, 2007. 

\item [{[21]}]    Rosen K. H., {\it et al}., {\it Handbook of Discrete and Combinatorial Mathematics}. CRC Press, 2000. 

\item [{[22]}]    Salas A. H., {\it About k-Fibonacci Numbers and their Associated Numbers}, International Mathematical Forum, {\bf 6}, (2011) No. 50, 2473-2479.  

\item [{[23]}]    Sigler L. E., {\it Fibonacci's Liber Abaci, A Translation into Modern English of Leonardo Pisano's Book of Calculation} Springer-Verlag New York, Inc., 2002. 

\item [{[24]}]    Spinadel V. W. de., {\it  The metallic means and design,} Nexus II: Architecture and Mathematics. Edizioni
dell’Erba, (1998), 141-157. 

\item [{[25]}]   Spinadel V. W. de., {\it  The family of metallic means,} Vis. Math., {\bf 1}, (1999), No. 3. 

\item [{[26]}] Spinadel V. W. de., {\it  The metallic means family and forbidden symmetries.} Int. Math. J, {\bf 2}, (2002), No. 3, 279-288. 

\end{itemize}

\end{document}